# The Most Beautiful Approach to the Most Beautiful Formula

## Amir Asghari

Asghari.amir@gmail.com

https://amirasghari.com/

The fact that $e^{i\pi} + 1 = 0$ is nothing less than magic. But, in the light of Euler's formula $e^{i\theta} = \cos\theta + i\sin\theta$, it is an easy exercise, taking $\theta = \pi$. However, according to the unwritten law of conversation of mathematical magic, the magic is not destroyed, it only moves to somewhere else, here, the Euler's justification of $e^{i\theta} = \cos\theta + i\sin\theta$. The trick is to put $i\theta$ in the Maclaurin series of $e^x$:

$$e^{i\theta} = 1 + i\theta + \frac{(i\theta)^2}{2!} + \frac{(i\theta)^3}{3!} + \frac{(i\theta)^4}{4!} + \frac{(i\theta)^5}{5!} + \frac{(i\theta)^6}{6!} + \cdots$$

Do the algebra and reorder as a complex number:

$$e^{i\theta} = (1 - \frac{\theta^2}{2!} + \frac{\theta^4}{4!} - \frac{\theta^6}{6!} + \cdots) + i(\theta - \frac{\theta^3}{3!} + \frac{\theta^5}{5!} - \cdots)$$

Voila! Here is our beautiful formula:

$$e^{i\theta} = \cos\theta + i\sin\theta$$

It is easy when we are guided to it. However, it has three undesirable aspects:

- it remains as an individual trick.
- It presupposes the formula.
- It unnecessarily links to the series (though, this aspect could be desirable when showing the beauty and usefulness of the series.)

It might be for these reasons that our beautiful formula couldn't find its rightful place in calculus textbooks, usually sitting somewhere in the appendices. This short note is to reclaim that rightful place while resolving all the aforementioned undesirable aspects.

**When the only trick is the beauty of connections**

When we multiply two powers with the same base, we add the exponents.

| $\theta_1$ | $\theta_2$ | $\theta_1 + \theta_2$ |
|---|---|---|
| $a^{\theta_1}$ | $a^{\theta_2}$ | $a^{\theta_1 + \theta_2}$ |

When we multiply two imaginary numbers with the modulus one, we add the arguments.

| $\theta_1$ | $\theta_2$ | $\theta_1 + \theta_2$ |
|---|---|---|
| $\cos\theta_1 + i\sin\theta_1$ | $\cos\theta_2 + i\sin\theta_2$ | $\cos(\theta_1 + \theta_2) + i\sin(\theta_1 + \theta_2)$ |

This suggests that the table of the unit imaginary numbers could be a table of powers. If so, what power is it?

$$cos\,\theta + i\sin\theta = a^b$$

We neither know *a* nor *b*. Here comes *the paradoxical essence of calculus*.

If you do not know what a function is, try to find how it changes!

This is what we do to find *a* and *b*, by taking the derivative of $cos\,\theta + i\sin\theta$.

$$(cos\,\theta + i\sin\theta)' = i(cos\,\theta + i\sin\varphi\theta)$$

That means:

$$(a^b)' = ia^b$$

And we know this works for $e^{i\theta}$ and Voila!

$$cos\,\theta + i\sin\theta = e^{i\theta}$$

(knowing that $\cos 0 + i\sin 0 = e^{i0}$)

Of course, we can make it more rigorous if we wish. For example, we can see both tables as examples of $f(x)f(y) = f(x + y)$, hence proving that the table of the unit imaginary numbers is a table of powers. But when teaching calculus, I favour beauty (of connections) over rigour.